\documentclass{llncs}
\usepackage[english]{babel}
\usepackage{amssymb,amsmath,latexsym} 
\usepackage[latin1]{inputenc}
\numberwithin{equation}{section}
\def\neu#1{{\bf #1}}
\def\tru{\triangle}
\def\trd{\bigtriangledown}
\def\tu{^\tru}
\def\td{^\trd}
\def\wdl{weakly dicomplemented lattice}

\def\WDL{\mbox{\rm WDL}}
\def\WDN{\mbox{\rm WDN}}
\def\WCL{\mbox{\rm WCL}}
\def\DCL{\mbox{\rm DCL}}
\def\Con{\mbox{\rm Con}}
\def\Id{\mbox{\rm Id}}

\title{Axiomatization of Boolean algebras via weak dicomplementations}
\author{L\'eonard Kwuida}
\institute{Universit\'e du Qu\'ebec en Outaouais\\ Gatineau, Canada\\
{\tt leonard.kwuida@uqo.ca}}

\begin{document}
\maketitle
\begin{abstract}
In this note we give an axiomatization of Boolean algebras based on weakly dicomplemented lattices: an algebra $(L,\wedge,\vee,\tu)$ of type $(2,2,1)$ is a Boolean algebra iff $(L,\wedge,\vee)$ is a non empty lattice and $(x\wedge y)\vee(x\wedge y\tu)=(x\vee y)\wedge(x\vee y\tu)$ for all $x,y\in L$. This provides a unique equation to encode distributivity and complementation on lattices. 
\end{abstract}
\section{Weakly dicomplemented lattices}\label{S:wdicomp}
\begin{definition}[\cite{Kw04}]\label{D:wdl}
A \neu{weakly dicomplemented lattice} is a bounded lattice $L$ equipped with two unary operations $\tu$ and $\td$ called \neu{weak complementation} and \neu{dual weak complementation}, and satisfying for all $x,y\in L$ the conditions:
\medskip\par\noindent
\begin{minipage}{.45\textwidth}
\begin{enumerate}
\item[(1)] $x^{\tru\tru}\le x$,
\item[(2)] $x\le y\implies x\tu\ge y\tu$,
\item[(3)] $(x\wedge y)\vee(x\wedge y\tu)=x$,
\end{enumerate}
\end{minipage}\hfill
\begin{minipage}{.45\textwidth}
\begin{enumerate}
\item[(1')] $x^{\trd\trd}\ge x$,
\item[(2')] $x\le y\implies x\td\ge y\td$,
\item[(3')] $(x\vee y)\wedge(x\vee y\td)=x$.
\end{enumerate}
\end{minipage}\medskip\par\noindent 
We call $x\tu$ the \neu{weak complement} of $x$ and $x\td$ the \neu{dual weak complement} of $x$. The pair $(x\tu,x\td)$ is called the \neu{weak dicomplement} of $x$ and the pair $(\tu,\td)$ a  \neu{weak dicomplementation} on $L$. The structure $(L,\wedge,\vee,\tu,0,1)$ is  called a \neu{weakly complemented lattice} and $(L,\wedge,\vee,\td,0,1)$ a \neu{dual weakly complemented lattice}.
\end{definition}
Note that $x^{\tru\tru}\le x\iff x^{\tru\tru}\vee x=x$ and $x^{\trd\trd}\ge x\iff x^{\trd\trd}\wedge x =x$; thus conditions (1) and (1') can be writen as equations. For (2) and (2'), we have $x\le y\implies x\tu\ge y\tu$ is equivalent to $(x\wedge y)\tu\wedge y\tu= y\tu$ and $x\le y\implies x\td\ge y\td$ equivalent to $(x\wedge y)\td\wedge y\td =y\td$. Therefore the class of weakly dicomplemented lattices form a variety. We denote it by \WDL. Similarly, the class \WCL\ of  weakly complemented lattices  and the class \DCL\ of dual weakly complemented lattices are varieties. These classes have been introduced to capture the notion of negation on ``concepts''\cite{Wi82,Wi96,Wi00,Kw04}, based on the work of Boole~\cite{Bo54}.

The following properties are easy to verify:
\begin{itemize}
\item[(4)] $y\vee y\tu=1$, $0\tu=1$,  $y\wedge y\td=0$, $1\td=0$ and $y\td\leq y\tu$ \qquad by (3) and (3'),
\item[(5)] $x\sp{\tru\tru\tru}=x\tu$ and $x\mapsto x\sp{\tru\tru}$ is a kernel operator on $L$ as well as $x\sp{\trd\trd\trd}=x\td$ and $x\mapsto x\sp{\trd\trd}$ is a closure operator on $L$ \qquad by (1)--(2'). 
\end{itemize}
Before we move to more properties, we give some examples.  
\begin{enumerate}
\item[(a)]The motivating examples of {\wdl}s are concept algebras. These are concept lattices with a weak negation and a weak opposition. For a detailed account on concept algebras, we refer the reader to \cite{Wi00}, \cite{Kw04} or \cite{GK07}.   
\item[(b)] The natural examples of {\wdl}s are Boolean algebras. In fact if $(B,\wedge,\vee,\bar{\phantom{a}},0,1)$ is a Boolean algebra then $(B,\wedge,\vee,\bar{\phantom{a}},\bar{\phantom{a}},0,1)$ (the complementation is duplicated, i.e. $x\tu := \bar{x}=: x\td$) is a weakly dicomplemented lattice.
\item[(c)]  Each bounded lattice can be endowed with a trivial weak dicomplementation by defining $(1,1),\,(0,0)$ and $(1,0)$ as the dicomplement of $0,\ 1$ and of each $x\not\in\{0,1\}$, respectively.
\end{enumerate}

\begin{theorem}\label{T:wcl}
Weakly complemented lattice are exactly non empty lattices satisfying the equations (1)--(3) in Definition~\ref{D:wdl}. 
\end{theorem}
Of course, weakly complemented lattices satisfy the equations (1)--(3) in Definition~\ref{D:wdl}. So what we should prove, is that, all lattices satifying the equations $(1)-(3)$ are bounded. 
\begin{proof}
Let $L$ be a non empty lattice satisfying the equations (1)--(3'). Let $x\in L$. We set $1:=x\vee x\tu$ and $0:=1\tu$. We are going to prove that $1$ and $0$ are respectively the greatest and lowest element of $L$. Let $y$ be an arbitrary element of $L$. We have $1\geq y\wedge 1=y\wedge(x\vee x\tu)\geq(y\wedge x)\vee(y\wedge x\tu)=y$, by (3). Thus $x\vee x\tu$ is the greatest element of $L$. Of course, if $L$ was equipped with a unary operation $\td$ satisfying the equation (1')--(3') we could use the same argument as above to say that $x\wedge x\td$ is the smallest element of $L$. Unfortunately we have to check that $0:=1\tu$ is less than every element of $L$. So let $y\in L$. We want to prove that $0\leq y$. Note that $(y\wedge y\tu)\tu\geq y\tu\vee y^{\tru\tru}=1$; thus $(y\wedge y\tu)\tu=1$. For an arbitrary element $z$ of $L$, we have 
\[0\wedge z=1\tu\wedge z=(y\wedge y\tu)^{\tru\tru}\wedge z\leq y\wedge y\tu\wedge z\leq y\wedge z\]
and
\[0\wedge z\tu=1\tu\wedge z\tu=(y\wedge y\tu)^{\tru\tru}\wedge z\tu\leq y\wedge y\tu\wedge z\tu\leq y\wedge z\tu.\]
Henceforth $0=(0\wedge z)\vee(0\wedge z\tu)\leq (y\wedge z)\vee(y\wedge z\tu)=y$. \qed
\end{proof}  
\def\thmref#1{Theorem~\ref{#1}}
\def\corref#1{Corollary~\ref{#1}}
\def\proref#1{Proposition~\ref{#1}}
\def\secref#1{Section~\ref{#1}}
\def\unl#1{\underline{#1}}
\def\Neu#1{\neu{#1}}

\section{Axiomatizing Boolean algebras}
\begin{definition}\label{D:wdl+neg}
 A weakly dicomplemented lattice is said to be \neu{with negation} if the unary operations coincide, i.e., if $x\td = x\tu$ for all $x$. 
   \end{definition}
The class of {\wdl}s with negation forms a subvariety of the class of all {\wdl}s, denoted by $\WDN$. 
\begin{theorem}
A {\wdl} $(L,\wedge,\vee,\tu,\td,0,1)$ is with negation iff $(L,\wedge,\vee,\tu,0,1)$ and $(L,\wedge,\vee,\td,0,1)$ are Boolean algebras. 
\end{theorem}
\begin{proof}
We assume that $(L,\wedge,\vee,\tu,\td,0,1)$ is a \wdl and $\tu=\td$. 
Recall that $x\vee x\tu=1$ and $x\wedge x\td=0$. Then with $\tu=\td$, $x\tu$ is a complement of $x$. To prove the distributivity, we will show that the lattices in the variety $\WDN$ of {\wdl}s with negation are all distributive. To this end is is enough to show that $\WDN$ is generated by the two element lattice, i.e  every member of $\WDN$ with at least three elements is not subdirectly irreducible.  We are going to show that for any $L\in\WDN$ with $|L|\geq 3$ there is $\theta_1,\theta_2\in\Con(L)$ such that $\theta_1\cap\theta_2=\Delta$, the trivial congruence (see for example \cite{BS81}).   
 \begin{itemize}
 \item[\textup(i)] For $c\in L\setminus\{0,1\}$, we have $[c,1]\cong[0,c\tu]$.  In fact the maps 
 \[
 \begin{array}{rcl}
 u_{c\tu} : [c,1] & \to & [0,c\tu]\\
  x & \mapsto & x\wedge c\tu
 \end{array} 
\qquad \text{ and }\qquad
 \begin{array}{rcl}
 v_{c} : [0,c\tu] & \to & [c,1]\\
  x & \mapsto & x\vee c
 \end{array} 
  \]
 are order preserving, injective and inverse of each other, since 
\[\left.
 \begin{array}{l}
x\wedge c\tu=y\wedge c\tu\\
x,y\geq c
 \end{array} 
\right\}\implies x=(x\wedge c)\vee(x\wedge c\tu)=c\vee(y\vee c\tu)=y,
 \]
\[\left.
 \begin{array}{l}
x\vee c=y\vee c\\
x,y\leq c\tu
 \end{array} 
\right\}\implies x=(x\vee c)\wedge(x\vee c\td)=(y\vee c)\wedge c\tu=y,
 \]
and  $v_c\circ u_{c\tu}(x)=(x\wedge c\tu)\vee c=(x\wedge c\tu)\vee(x\wedge c)=x=\Id_{[c,1]}(x)$ as well as $u_{c\tu}\circ v_c(x)=(x\vee c)\wedge c\tu=(x\vee c)\wedge(x\vee c\td)=x=\Id_{[0,c\tu]}(x)$.
 \item[\textup(ii)] The maps 
\[
\begin{array}{rcrcl}
 f_1 : L & \to & [c\tu,1] & \to & [0,c^{\tru\tru}]=[0,c]\\
  x & \mapsto & x\vee c\tu & \mapsto & (x\vee c\tu)\wedge c=x\wedge c
 \end{array} 
\]
and 
\[
 \begin{array}{rcrcl}
 f_2 : L & \to & [c,1] & \to & [0,c\tu]\\
  x & \mapsto & x\vee c & \mapsto & (x\vee c)\wedge c\tu=x\wedge c\tu
 \end{array} 
 \]
are lattice homomorphisms. In fact $f_1$ and $f_2$ trivially preserve $\wedge$; For $x,y$ in $L$ we have, $f_1(x)=u_{c\sp{\tru\tru}}(x\vee c\tu)$ and $f_2(x)=u_{c\tu}(x\vee c)$. In addition, 
\begin{eqnarray*}
f_1(x\vee y)&=& u_{c\sp{\tru\tru}}(x\vee y\vee c\tu)=u_{c\sp{\tru\tru}}\left((x\vee c\tu)\vee(y\vee c\tu)\right)\\
& = & u_{c\sp{\tru\tru}}(x\vee c\tu)\vee u_{c\sp{\tru\tru}}(y\vee c\tu)\\
&=& f_1(x)\vee f_1(y),
\end{eqnarray*}
and 
\begin{eqnarray*}
f_2(x\vee y)&=& u_{c\tu}(x\vee y\vee c)=u_{c\sp{\tru}}\left((x\vee c)\vee(y\vee c)\right)\\
& = & u_{c\sp{\tru}}(x\vee c)\vee u_{c\sp{\tru}}(y\vee c)\\
&=& f_2(x)\vee f_2(y). 
\end{eqnarray*}

\item[\textup(iii)] We set $\theta_1:={\ker}f_1$ and $\theta_2:={\ker}f_2$.
Then $\theta_1\cap\theta_2=\Delta$. In fact 
\begin{eqnarray*}
(x,y)\in\theta_1\cap\theta_2 &\implies & x\wedge c=y\wedge c \text{ and } x\wedge c\tu=y\wedge c\tu\\
&\implies& x=(x\wedge c)\vee(x\wedge c\tu)=(y\wedge c)\vee(y\wedge c\tu)=y\\
&\implies& (x,y)\in\Delta.
\end{eqnarray*}
 \end{itemize}\qed
\end{proof}

 \begin{corollary}\label{C:C1}
 $(L,\wedge,\vee,\tu,0,1)$ is a Boolean algebra iff $(L,\wedge,\vee)$ is a non empty lattice in which $x^{\tru\tru}=x$, $x\leq y\implies x\tu\geq y\tu$ and $(x\wedge y)\vee(x\wedge y\tu)=x=(x\vee y)\wedge(x\vee y\tu)$. 
 \end{corollary}

 \begin{theorem}[New axiom for Boolean algebras]\label{T:main}
 An algebra $(L,\wedge,\vee,\tu,0,1)$ is a Boolean algebra iff $(L,\wedge,\vee)$ is a non empty lattice in which 
\[(x\wedge y)\vee(x\wedge y\tu)=(x\vee y)\wedge(x\vee y\tu) \text{ for all }x,y\in L \qquad \textup{(\ddag)}.
\]
 \end{theorem}
 \begin{proof}
We are going to show that the equations in Corollary~\ref{C:C1} can be derived from (\ddag). 
\begin{itemize}
 \item[\textup(i)] $x\geq(x\wedge y)\vee(x\wedge y\tu)=(x\vee y)\wedge(x\vee y\tu)\geq x$ implies $(x\wedge y)\vee(x\wedge y\tu)=x=(x\vee y)\wedge(x\vee y\tu)$. 
 \item[\textup(ii)] $x=(x\vee y)\wedge(x\vee y\tu)$ implies $y\wedge y\tu=0$; thus $x=(x\wedge x\tu)\vee(x\wedge x^{\tru\tru})=0\vee(x\wedge x^{\tru\tru})=x\wedge x^{\tru\tru}$. Hence $x\leq x\sp{\tru\tru}$.
 \item[] $x=(x\wedge y)\vee(x\wedge y\tu)$ implies $y\vee y\tu=1$; thus $x=(x\vee x\tu)\wedge(x\vee x^{\tru\tru})=1\wedge(x\vee x^{\tru\tru})=x\vee x^{\tru\tru}$. Hence $x\geq x^{\tru\tru}$. Therefore $x=x^{\tru\tru}$.
 \item[\textup(iii)] Let $x\leq y$. Then $x\vee x\tu=1$ implies $y\vee x\tu=1$. Thus $x\tu=(x\tu\vee y\tu)\wedge(x\tu\vee y^{\tru\tru})=(x\tu\vee y\tu)\wedge(x\tu\vee y)=x\tu\vee y\tu$, and $x\tu\geq y\tu$. 
 \end{itemize}
\qed
 \end{proof}
In the proof of Theorem~\ref{T:main}, we have shown that the conditions (1)--(2') in Definition~\ref{D:wdl} follow from (3) and (3'), in case $\tu=\td$. Does this hold in general?


\begin{thebibliography}{00}
\bibitem[Bo54]{Bo54} G.\,Boole. \newblock \textsl{An investigation of the laws of thought on which are founded the ma\-the\-matical theories of logic and probabilities}. \newblock Macmillan 1854. Reprinted by Dover Publ. New york (1958).
\bibitem[BS81]{BS81} S.\,Burris \& H.\,P.\,Sankappanavar. \newblock \textsl{A course in universal algebra}. \newblock Springer Verlag (1981).
\bibitem[GK07]{GK07} B.\,Ganter \& L.\,Kwuida. \newblock \textsl{Finite distributive concept algebras}. \newblock Order. (2007).
\bibitem[Kw04]{Kw04} L.\,Kwuida. \newblock \textsl{Dicomplemented lattices. A contextual generalization of Boolean algebras}. \newblock Dissertation TU Dresden. Shaker Verlag. (2004).
\bibitem[Wi82]{Wi82} R.\,Wille. \newblock \textsl{Restructuring lattice theory:  an approach based on hierarchies of concepts}. \newblock in I. Rival (Ed.) \newblock \textsl{Ordered Sets}. Reidel (1982) 445-470.
\bibitem[Wi96]{Wi96} R.\,Wille. \newblock \textsl{Restructuring mathematical logic: an approach based on Peirce's pragmatism}. in  Ursini, Aldo (Ed.) \newblock \textsl{Logic and algebra} Marcel Dekker. Lect. Notes Pure Appl. Math. \textbf{180} (1996) 267-281.
\bibitem[Wi00]{Wi00} R.\,Wille.  \newblock \textsl{Boolean Concept Logic} in B. Ganter \& G.W. Mineau (Eds.) ICCS 2000 \newblock \textsl{Conceptual Structures: Logical, Linguistic, and Computational Issues} Springer LNAI \textbf{1867} (2000) 317-331.
\end{thebibliography}
\end{document}